\documentclass[11pt]{amsart}

\usepackage{amsmath, amssymb}
\usepackage{xr}

\DeclareMathOperator{\dist}{dist}
\newcommand{\fc}{\mathfrak c}

\newcommand{\bbN}{{\mathbb N}}

\newcommand{\bbR}{\mathbb R}
\newcommand{\bbC}{\mathbb C}

\newcommand{\cK}{{\mathcal K}}

\DeclareMathOperator{\id}{id}

\newcommand{\cC}{\mathcal C}

\newcommand{\bbP}{\mathbb P}
\newcommand{\e}{\varepsilon}
\newtheorem{theorem}{Theorem}[section] 

\newtheorem{lemma}[theorem]{Lemma}
\newtheorem{prop}[theorem]{Proposition}

\theoremstyle{definition}

\newtheorem{definition}[theorem]{Definition}

\newtheorem{problem}[theorem]{Problem}

\newcounter{my_enumerate_counter}
\newcommand{\pushcounter}{\setcounter{my_enumerate_counter}{\value{enumi}}}
\newcommand{\popcounter}{\setcounter{enumi}{\value{my_enumerate_counter}}}

\newcommand{\lbl}{\label}

\newcommand{\cP}{{\mathcal P}}
\newcommand{\cO}{{\mathcal O}}

\newcommand{\bbZ}{\mathbb Z}

\DeclareMathOperator{\Ad}{Ad}

\DeclareMathOperator{\Aut}{Aut}

\newcommand{\sn}{\kappa}
\newcommand{\sntoo}{\kappa'}

\DeclareMathOperator{\Cu}{Cu}
\DeclareMathOperator{\INV}{INV}

\title{Nonseparable UHF algebras II: Classification}

\author{Ilijas Farah}

\address{Department of Mathematics and Statistics\\
York University\\
4700 Keele Street\\
North York, Ontario\\ Canada, M3J 1P3\\
and Matematicki Institut, Kneza Mihaila 35, Belgrade, Serbia}

\urladdr{http://www.math.yorku.ca/$\sim$ifarah}

\email{ifarah@mathstat.yorku.ca}

\author{Takeshi Katsura}

\address{Department of Mathematics\\
Faculty of Science and Technology\\
Keio University\\
3-14-1 Hiyoshi, Kouhoku-ku, Yokohama\\
JAPAN, 223-8522}

\email{katsura@math.keio.ac.jp}

\date{\today.}
\begin{document}
\begin{abstract}
For every uncountable cardinal $\kappa$ there are $2^\kappa$ 
nonisomorphic  simple AF algebras of density character $\kappa$ and 
$2^\kappa$ nonisomorphic hyperfinite II$_1$ factors of density character $\kappa$.
These estimates are maximal possible. 
All C*-algebras that we construct 
 have the same Elliott invariant and Cuntz semigroup as the CAR algebra.    \end{abstract}
\maketitle

\section{Introduction}

The classification program of nuclear separable C*-algebras 
can be traced back to classification of UHF algebras of Glimm 
and Dixmier. However, it was Elliott's classification of 
AF algebras  and real rank zero  AT algebras that started the classification 
program in earnest (see  e.g., \cite{Ror:Classification} and \cite{EllTo:Regularity}).

While it was generally agreed that the
classification of nonseparable C*-algebras is a nontractable problem, there
were no concrete results to this effect. 
Methods from logic were recently successfully applied to analyze the 
classification problem for separable C*-algebras (\cite{FaToTo:Turbulence}) 
and II$_1$ factors with separable predual (\cite{sato09b}) and it comes as no 
surprise that they are also instrumental in analyzing classification of nonseparable
operator algebras.  
 We construct large families of
nonseparable AF algebras with identical K-theory and Cuntz semigroup as the
CAR algebra. 
 Since the CAR algebra is a prototypical example of a classifiable algebra, 
this gives a strong endorsement  to the above viewpoint. 
We also construct a large family of hyperfinite II$_1$ factors with predual 
of character density $\kappa$ 
for every uncountable cardinal $\kappa$. 
Recall that a \emph{density character} of a metric space is the least
 cardinality of a dense subset.  While the CAR algebra is unique and 
 there is a unique hyperfinite II$_1$ factor with 
 separable predual, our results show that  uniqueness badly fails 
 in every uncountable density character 
 $\kappa$.

For each $n
\in \bbN$, we denote by $M_n(\bbC)$ the unital C*-algebra of all
$n\times n$ matrices with complex entries. A C*-algebra which is
isomorphic to $M_n(\bbC)$ for some $n \in \bbN$ is called a
\emph{full matrix algebra}.

\begin{definition}
A C*-algebra $A$ is said to be
\begin{itemize}
\item \emph{uniformly hyperfinite}
(or \emph{UHF})
if $A$ is isomorphic to a tensor product
of full matrix algebras.
\item \emph{approximately matricial}
(or \emph{AM})
if it has a directed family of
full matrix subalgebras with dense union.
\item locally matricial (or LM) if for any finite subset $F$ of $A$ and any $\e > 0$, there exists a full
matrix subalgebra $M$ of $A$ such that $\dist(a,M)<\e$ for all $a\in F$. 
\end{itemize}
\end{definition}
  In \cite{Dix:Some}
Dixmier remarked that in the unital case 
these three classes coincide under the additional assumption that $A$ is separable 
and asked whether this result extends to nonseparable algebras. 
In   \cite{Kat:Non-separable} a pair of nonseparable AF algebras not isomorphic
to each other but with the same Bratteli diagram was constructed.
Dixmier's question was answered in the negative in 
 \cite{FaKa:Nonseparable}.  Soon after, 
 AM algebras with counterintuitive properties were constructed.  
A simple nuclear algebra that has irreducible representations 
on both separable and nonseparable Hilbert space 
was constructed in   \cite{Fa:Graphs} and an algebra 
with nuclear dimension zero which does not absorb the Jiang--Su algebra tensorially was constructed in \cite{FaHaTiKa:Simple}. Curiously, all of these 
results (with the possible exception of   \cite{FaHaTiKa:Simple})
were proved in ZFC.

 Results of the  present paper widen the gap between unital UHF and 
 AM algebras even further  by showing that there are many more AM algebras than UHF algebras of  every uncountable  density character. 
 In \S\ref{Sec:Non-classAM-1}~and~\S\ref{Sec:Non-classAM-2}
 we prove the following.

\begin{theorem}\label{T.non-classification}
For every uncountable cardinal $\kappa$ 
 there are $2^\kappa$ pairwise nonisomorphic AM algebras with character
density $\kappa$. All these algebras have the same $K_0$, $K_1$, and Cuntz
semigroup as the CAR algebra.\end{theorem}

Every AM algebra is LM and by Theorem~\ref{T.non-classification} 
 there are already as many AM algebras as there are C*-algebras in every uncountable density character. Therefore 
 no quantitive information along these lines can be obtained about LM algebras.

\begin{theorem}\label{T.non-classification.hyperfinite} 
For every uncountable cardinal $\kappa$ there are $2^\kappa$
nonisomorphic hyperfinite II$_1$ factors with predual of density character $\kappa$. 
\end{theorem}

While there is a unique hyperfinite II$_1$ factor with separable predual, 
it was proved by Widom (\cite{Wid:Nonisomorphic}) that there are at least
as many nonisomorphic hyperfinite II$_1$ factors with predual of density character $\kappa$
as there are infinite cardinals $\leq\kappa$.

Note that there are at most $2^\kappa$
C*-algebras of density character $\kappa$ and at most $2^\kappa$ von Neumann 
algebras with predual of density character $\kappa$.
This is because each such algebra has a dense  subalgebra of cardinality $\kappa$, 
and an easy counting argument shows that there
are at most $2^\kappa$ ways to define $+, \cdot, {}^*$ and $\|\cdot\|$ on a
fixed set of size $\kappa$.

On the positive side, in Proposition~\ref{P.UHF-classification} we show that 
Glimm's classification of UHF algebras by their generalized integers 
extends to nonseparable algebras. This shows that the number of 
isomorphism classes of UHF algebras of density character $\leq \kappa$
is equal to  $2^{\aleph_0}$, as long as there are only countably many cardinals $\leq \kappa$
(Proposition~\ref{P.UHF.counting} and the table in \S\ref{Sec:Open}). 
Hence UHF algebras of arbitrary 
density character  are `classifiable' in the sense of Shelah (e.g., \cite{Sh:Classification}). 
Note, however, that they don't form an elementary class (cf. \cite{Mitacs2012}). 

Two  C*-algebras are isomorphic if and only if they are isometric, and the same fact is 
true for  II$_1$ factors with $\ell_2$-metric. However, in some situations there exist topologically isomorphic but not
isometric structures---notably, in the case of Banach spaces.    
The more general problem of constructing many nonisomorphic models in a given 
density character was considered in~\cite{ShUs:928}.

\subsection*{Organization of the paper}
In \S\ref{Sec:Pre} we set up the toolbox used in the paper.  
In \S\ref{Sec:K}
we study K-theory and Cuntz semigroup of nonseparable LM algebras. 
UHF algebras are classified in \S\ref{Sec:Classification}. 
In 
\S\ref{Sec:Non-classAM-1} we prove a non-classification result for 
AM algebras and hyperfinite II$_1$ factors in regular character densities. 
Shelah's methods from \cite{Sh:E59}, as adapted to the context of metric structures
in~\cite{FaSh:Dichotomy} 
are used to extend this to arbitrary uncountable character densities
in~\S\ref{Sec:Non-classAM-2}.
In \S\ref{Sec:Open} we state some open problems and provide some limiting examples. 

The paper requires only basic background in operator algebras (e.g., \cite{Black:Operator}) 
and in naive set theory. On several occasions we include remarks aimed at 
model theorists. Although they provide an 
additional insight, these  remarks can be safely ignored by readers not interested in model theory. 

\subsection*{Acknowledgments} 
Results of the present paper were proved at the Fields Institute in January  2008 (the case when $\kappa$ is a regular cardinal) and at the Kyoto University in November 2009. We would like to thank 
both institutions for their hospitality. 
I.F. would like to thank to Aaron Tikuisis for many useful remarks on the draft of this paper 
and  for correcting the proof of Proposition~\ref{P.5.7},  and to  
 Teruyuki Yorioka for supporting his visit to Japan. 
We would also like to thank David Sherman for providing reference to Widom's paper \cite{Wid:Nonisomorphic}. 
I.F. is partially supported by NSERC.

\section{Preliminaries}\label{Sec:Pre}

A cardinal $\kappa$ is a \emph{successor cardinal} if it is the least cardinal
greater than some other cardinal. A cardinal that is not a successor is called
a \emph{limit cardinal}. Note that every infinite cardinal is a limit ordinal.
  Cardinal~$\kappa$ is
\emph{regular} if for $X\subseteq \kappa$ we have $\sup X=\kappa$ if and only
if $|X|=\kappa$. For example, every successor cardinal is regular. 
A cardinal that is not regular is \emph{singular}. The least singular cardinal 
is $\aleph_\omega$ and singular cardinal combinatorics is a notoriously difficult subject. 
A subset $C$
of an ordinal $\gamma$ is \emph{closed and unbounded} (or \emph{club}) if its
supremum is $\gamma$ and whenever  $\delta<\gamma$ is such that $\sup(C\cap
\delta)=\delta$ we have $\delta\in C$. A subset of an ordinal  $\gamma$  is
called stationary if it intersects every club in $\gamma$ non-trivially.

Some of the lemmas in the present paper, 
 (e.g., Lemma~\ref{L.many-subsets}) are well-known  
but we provide proofs for the convenience of the readers.

\begin{lemma} \label{L.many-subsets}
If $\kappa$ is a regular cardinal then there
are $S(X)\subseteq \kappa$, $X\subseteq \kappa$ such that the symmetric
difference $S(X)\Delta S(Y)$ is stationary whenever $X\neq Y$.
\end{lemma}

\begin{proof}
We first prove that $\kappa$ can be partitioned into $\kappa$ many stationary
sets, $Z_\gamma$, $\gamma<\kappa$.  If $\kappa$ is a successor cardinal then
this is a result of Ulam (\cite[Corollary~6.12]{Ku:Book}). If $\kappa$ is a
limit cardinal, then there are $\kappa$ regular cardinals below $\kappa$. For
each such cardinal the set  $Z_\gamma=\{\delta<\kappa : \min\{|X|: X\subseteq
\delta$ and $\sup X=\delta\}=\gamma\}$ is stationary.

 For $X\subseteq
\kappa$ let $S(X)=\bigcup_{\gamma\in X} Z_\gamma$. Then clearly $S(X)\Delta
S(Y)$ is stationary whenever $X\neq Y$.
\end{proof}

Let $|X|$ denote the cardinality of a set $X$. 
We shall now recall some basic set-theoretic notions 
worked out explicitly in the case of C*-algebras in \cite{FaKa:Nonseparable}. 

\begin{definition}
A directed set $\Lambda$ is said to be \emph{$\sigma$-complete}
if every countable directed $Z\subseteq \Lambda$
has the supremum $\sup Z \in \Lambda$.
A directed family $\{A_\lambda\}_{\lambda\in \Lambda}$
of subalgebras of a C*-algebra $A$
is said to be \emph{$\sigma$-complete}
if $\Lambda$ is $\sigma$-complete
and for every countable directed $Z\subseteq \Lambda$,
$A_{\sup Z}$ is the closure of the union of $\{A_\lambda\}_{\lambda\in Z}$.
\end{definition}

Assume $A$ is a  nonseparable C*-algebra.  Then 
$A$ is a direct limit of a $\sigma$-complete directed system of
its separable subalgebras (\cite[Lemma~2.10]{FaKa:Nonseparable}). 
Also, if $A$ is represented as a direct limit of a $\sigma$-complete directed 
system of separable subalgebras in two different ways, then 
the intersection of these two systems is a $\sigma$-complete directed system 
of separable subalgebras and~$A$ is its direct limit (\cite[Lemma~2.6]{FaKa:Nonseparable}). 

The following was proved in \cite[remark following Lemma~2.13]{FaKa:Nonseparable}. 

\begin{lemma}\lbl{L.AM-reflection}
A C*-algebra $A$ is LM if and only if it is equal to a union of 
a $\sigma$-complete directed family of separable AM subalgebras. \qed
\end{lemma}

In \S\ref{Sec:Non-classAM-2}
we shall use the following well-known fact without mentioning.
We give its proof for the reader's convenience.

\begin{lemma}
Let $\alpha$ be an action of a group $G$ on a unital C*-algebra $A$. Let
$\{u_g\}_{g\in G} \subset A \rtimes_\alpha^r G$ be the implementing unitaries in the reduced crossed product.
Suppose that
a unital subalgebra $A_0 \subset A$ and a subgroup $G_0 \subset G$ satisfy that
$\alpha_g[A_0] = A_0$ for all $g \in G_0$,
and set $B_0 := C^*(A_0 \cup \{u_g\}_{g \in G_0})$.
Then we have
\[
B_0 \cap A = A_0 \qquad\text{ and }\qquad 
B_0 \cap \{u_g\}_{g \in G} = \{u_g\}_{g \in G_0}
\]
in $A \rtimes_\alpha G$.
\end{lemma}

\begin{proof}
First note that there exists a conditional expectation $E$ onto $A \subset A
\rtimes_\alpha G$ such that $E(a)=a$ and $E(a u_g) = 0$ for all $a \in A$ and $g \in G
\setminus \{e\}$  (see \cite[Proposition~4.1.9]{BrOz}). Since the linear span of $\{a u_g : a \in
A_0, g \in G_0\}$ is dense in $B_0$, we have $E[B_0] = A_0$. This shows $B_0
\cap A = E[B_0 \cap A] =A_0$. For the same reason we have $E[B_0u_g^*] = 0$ for
all $g \in G \setminus G_0$. This shows that $u_g \notin B_0$ for $g \in G
\setminus G_0$. Thus $B_0 \cap \{u_g\}_{g \in G} = \{u_g\}_{g \in G_0}$.
\end{proof}

\section{$K$-theory of LM algebras}
\label{Sec:K}

For definition of groups $K_0(A)$ and $K_1(A)$  see e.g.,
\cite{Black:Operator} or \cite{RoLaLa:Introduction}
and for the Cuntz semigroup $\Cu(A)$ see e.g., \cite{CoElIv}. 
 If $A$ is a unital subalgebra of $B$ then $K_1(A)$ is a subgroup
of $K_1(B)$ and if $B=\varinjlim B_\lambda$ then $K_1(B)=\varinjlim
K_1(B_\lambda)$. Both of these two properties fail for $K_0$ in general.

A reader familiar with the logic of metric structures (\cite{BYBHU}, \cite{FaHaSh:Model1}) 
will notice that in  Lemma~\ref{L.K-0.reflection} we are only using two standard facts: 
(1) the family of separable elementary submodels of algebra $A$ is $\sigma$-complete
and has $A$ as its direct limit and (2)  
if $A_\lambda$ is an elementary submodel of $A$ then $K_0(A_\lambda)$ is a subgroup of $K_0(A)$
and $\Cu(A_\lambda)$ is a subsemigroup of $\Cu(A)$.  

\begin{lemma} 
\label{L.K-0.reflection}
If $A$ is a nonseparable C*-algebra then $A$ is
a union of  a $\sigma$-complete directed family of separable subalgebras
$A_\lambda$, $\lambda\in \Lambda$, such that for each $\lambda\in \Lambda$ we
have
\begin{enumerate}
\item $K_0(A_\lambda)$ is a subgroup of $K_0(A)$ and $K_0(A)=\varinjlim
K_0(A_\lambda)$,
\item  $\Cu(A_\lambda)$ is a sub-semigroup of
$\Cu(A)$ and $\Cu(A)=\varinjlim \Cu(A_\lambda)$.
\end{enumerate}
\end{lemma}

\begin{proof} (1)  
As usual $p\sim q$ denotes the Murray--von Neumann equivalence of projections in algebra $A$, 
namely $p\sim q$ if and only if $p=vv^*$ and $q=v^*v$ for some $v$ in $A$. 

For a subalgebra $B$ of $A$
we have that $K_0(B)<K_0(A)$ if and only if for any two projections $p$ and $q$
in $B\otimes \cK$ we have $p\sim q$ in $B$ if and only if $p\sim q$ in $A$.

We need to show that the family of separable subalgebras $B$ of $A$ such that
$K_0(B)<K_0(A)$ is closed and unbounded. Since $\|p-q\|<1$ implies $p\sim q$,
this set is closed. The following condition for all $p,q$ in $B$ implies
$K_0(B)$ is a subgroup of $K_0(A)$:
$$
\inf_{v\in B} \|vv^*-p\|+\|v^*v-q\| =\inf_{v\in A} \|vv^*-p\|+\|v^*v-q\|.
$$
Since $\|p-p'\|<1$ implies $p\sim p'$,  a closing-up argument like in the proof
of \cite[Lemma~\ref{UHF1-L.AF-reflection}]{FaKa:Nonseparable}
 shows every separable subalgebra of $A$ is
contained in one that satisfies the above condition for each pair of
projections.

The assertion that $K_0(A)=\varinjlim K_0(A_\lambda)$ is automatic since
$A=\bigcup_\lambda A_\lambda$.

(2) Recall that the Cuntz ordering on positive elements in algebra $A$ is defined
by $a\precsim b$ if for every $\e>0$ there exists $x\in A$ such that $\|a-xbx^*\|<\e$. 

We need to show that the family of separable subalgebras $B$ of $A$ such
that for all $a$ and $b$ in $B$ we have $a\precsim b$ in $B$ if and only if
$a\precsim b$ in $A$ is closed and unbounded. It is clearly closed. Again it
suffices to assure that for a dense set of pairs $a,b$ of positive operators in
$B$ we have $\inf_{x\in B} \|a-xbx^*\|=\inf_{x\in A} \|a-xbx^*\|$, and this is
achieved by a L\"owenheim--Skolem argument resembling one in the proof of
Lemma~\ref{L.AM-reflection}.

The assertion that $\Cu(A)=\bigcup_\lambda \Cu(A_\lambda)$ is again automatic.
\end{proof}





It is also true that if $A$ is a nonseparable C*-algebra with the unique trace
then its separable subalgebras with the unique trace form a $\sigma$-complete
directed system whose direct limit is equal to $A$. This follows
from an argument due to N.C. Phillips (see \cite{Phi:Simple}) and it can be
proved by the argument of Lemma~\ref{L.K-0.reflection} 
(see also
\cite[Remark~\ref{UHF1-R.K-0.reflection}]{FaKa:Nonseparable}).

 Recall that $n$ is a
\emph{generalized integer} (or a supernatural number) 
if $n=\prod_{p\in \cP} p^{n_p}$ where $n_p\in
\bbN\cup\{\infty\}$ for all $p$.  For a unital UHF algebra $A$ define the
generalized integer $n=\prod_{p\in \cP}p^{n_p}$ of $A$ by
\[
n_p := \sup\big\{k\in \bbN: \text{there exists a unital homomorphism from
$M_{p^k}(\bbC)$ to $A$}\big\}
\]
for each $p \in \mathcal{P}$.

Glimm (\cite{Glimm:On}) has shown that the generalized integer provides a
complete invariant for isomorphism of separable unital UHF algebras. For a
generalized integer $n$ define the group
$$
\bbZ[1/n]=\{k/m: k\in \bbZ, m\in \bbZ\setminus \{0\}, m|n\}
$$
where $m|n$ is defined in the natural way. Then for a separable UHF algebra $A$
and its generalized integer $n$ we have $K_0(A)=\bbZ[1/n]$.

\begin{prop}\label{P.tracial-state}
An LM algebra $A$ has a unique tracial state $\tau$. If $A$ is unital, then
$\tau$ induces an isomorphism from $K_0(A)$ onto $\bbZ[1/n] \subset \bbR$, with
$n$ defined as above,  as ordered groups. We have $K_1(A) = 0$.
\end{prop}

\begin{proof}
Uniqueness 
of the tracial state immediately follows from the fact that a nonseparable LM algebra is a $\sigma$-complete direct limit of separable UHF algebras, since they have a unique tracial state.   
If $A$ is unital we fix $\tau$ so that
$\tau(1)=1$.

 For projections $p$ and $q$ of $A$ we have $\tau(p)=\tau(q)$
 if and only if $p\sim q$. This is true for separable LM algebras and the
 nonseparable case follows immediately by Lemma~\ref{L.AM-reflection}. Therefore $\tau$ is an isomorphic embedding of $K_0(A)$ into
$\bbZ[1/n]$. Since $K_1(B)=0$ for each separable LM algebra $A=\varinjlim
A_\lambda$ implies $K_1(A)=\varinjlim K_1(A_\lambda)$, we have $K_1(A)=0$ by
Lemma~\ref{L.AM-reflection}.
\end{proof}

The following is an immediate consequence of the main result of \cite{BrPeTo:Cuntz}. 

\begin{prop} 
If  $A$ is an infinite-dimensional LM algebra 
then  its Cuntz semigroup  is isomorphic to $K_0(A)_+\sqcup (0,\infty)$. \qed
\end{prop}

\section{Classification of UHF algebras}
\label{Sec:Classification}

\begin{lemma}\label{L.tensor.isomorphism}
Assume $A=\bigotimes_{x\in X} A_x$, $B=\bigotimes_{y\in Y} B_y$
and all $A_x$ and all $B_y$ are unital, separable, simple, and not equal to $\bbC$.
Let $\Phi\colon A\to B$ be an isomorphism.
Then there exist partitions $X = \bigsqcup_{z \in Z}X_z$ and $Y = \bigsqcup_{z \in
Z}Y_z$ of $X$ and $Y$ into disjoint nonempty countable subsets
indexed by a same set $Z$ such that
\[
\Phi[\textstyle\bigotimes_{x\in X_z} A_x]=\bigotimes_{y\in Y_z} B_y
\]
for all $z \in Z$.
\end{lemma}

\begin{proof}
Consider the set $\bbP$ of pairs of families $(\{X_z\}_{z \in Z},\{Y_z\}_{z \in Z})$
of disjoint nonempty countable subsets of $X$ and $Y$, respectively,
indexed by a same set $Z$
such that we have
$\Phi[\bigotimes_{x\in X_z}A_x]=\bigotimes_{y\in Y_z} B_y$
for every $z\in Z$.
Order $\bbP$ by letting
\[
(\{X_z\}_{z \in Z},\{Y_z\}_{z \in Z})
\leq 
(\{X_z'\}_{z \in Z'},\{Y_z'\}_{z \in Z'})
\]
if $Z\subseteq Z'$ and $X_z'=X_z$ and $Y_z'=Y_z$ for all $z\in Z$. 

By Zorn's lemma,
there exists a maximal one
$\{X_z\}_{z \in Z}$ and $\{Y_z\}_{z \in Z}$
among such families.
If we set $X':=X\setminus \bigcup_{z\in Z} X_z$
and $Y':=Y\setminus \bigcup_{z\in Z} Y_z$
then 
$\bigotimes_{x\in X'}A_x=\bigcap_{z\notin X'} Z_A(A_z)$
and 
$\bigotimes_{y\in Y'}B_y=\bigcap_{z\notin Y'} Z_B(B_z)$
by \cite[Theorem~1]{HaWa}. 
Therefore $\Phi[\bigotimes_{x\in X'}A_x]
=\bigotimes_{y\in Y'} Y_z$.
Thus $X'$ is nonempty if and only if $Y'$ is nonempty.
Suppose, to derive a contradiction,
both $X'$ and $Y'$ are nonempty.
By applying the argument in
the proof of \cite[Lemma~\ref{UHF1-Lem:phi=id}]{FaKa:Nonseparable}
 (see also \cite[Lemma~\ref{UHF1-Lem:tensor}]{FaKa:Nonseparable}),
we find non-empty
countable $X_0\subseteq X'$ and $Y_0\subseteq Y'$ such that
$\Phi[\bigotimes_{x\in X_0} A_x]=\bigotimes_{y\in Y_0} B_y$.
This contradicts the assumed maximality of
$\{X_z\}_{z \in Z}$ and $\{Y_z\}_{z \in Z}$.
Hence both $X'$ and $Y'$ are empty,
and the maximal families
$\{X_z\}_{z \in Z}$ and $\{Y_z\}_{z \in Z}$
are what we want.
\end{proof}

 Let us denote  the set of all prime
numbers by $\mathcal{P}$.

\begin{prop}\label{P.UHF-classification}
If $\kappa_p$, $\lambda_p$, $p\in \cP$ are sequences of cardinals indexed by
the prime numbers then $\bigotimes_{p\in \cP}\bigotimes_{\kappa_p} M_p(\bbC)$ and
$\bigotimes_{p\in \cP}\bigotimes_{\lambda_p} M_p(\bbC)$ are isomorphic if and only if
$\kappa_p=\lambda_p$ for all $p$.
\end{prop}

\begin{proof} Only the direct implication requires a proof.
The separable case is a theorem of Glimm (\cite{Glimm:On}). Assume the algebras
are nonseparable, and let $X=\bigcup_{p\in \cP}\{p\}\times \kappa_p$,
$A_{(p,\gamma)}=M_p(\bbC)$, $Y=\bigcup_{p\in \cP} \{p\}\times \lambda_p$, and
$B_{(p,\gamma)}=M_p(\bbC)$. By 
Lemma~\ref{L.tensor.isomorphism} applied to the
isomorphism between $\bigotimes_{x\in X} A_x$ and $\bigotimes_{y\in Y} B_y$ we
can find partitions $X=\bigcup_{z\in Z} X_z$ and $Y=\bigcup_{z\in Z} Y_z$ into
countable sets such that $\bigotimes_{x\in X_z} A_x$ and $\bigotimes_{y\in Y_z}
B_y$ are isomorphic for each $z\in Z$. By Glimm's theorem and simple cardinal
arithmetic this implies $\kappa_p=\lambda_p$ for all~$p$.
\end{proof}

By Proposition~\ref{P.UHF-classification}, for each UHF algebra
$A=\bigotimes_{p\in \cP}\bigotimes_{\kappa_p} M_p(\bbC)$ we can define the
generalized integer $\sn(A)=\prod_{p\in \cP} p^{\kappa_p}$ and UHF algebras are
completely classified up to isomorphism by the generalized integers $\sn(A)$
associated with them. Note that $\sn(A)$ being well-defined hinges on
Proposition~\ref{P.UHF-classification}. It is unclear whether $\sn(A)$
coincides with the generalized integer obtained by a straightforward
generalization of definition given for separable UHF algebras before
Proposition~\ref{P.tracial-state}; see Problem~\ref{Prob.sn} and
Problem~\ref{Prob.emb}. We shall avoid using this notation for generalized 
integers in order to avoid the confusion with powers of cardinal numbers.

\begin{prop}\label{P.UHF.counting}
For every ordinal $\gamma$ 
there are $(|\gamma|+\aleph_0)^{\aleph_0}$ isomorphism classes of unital 
UHF algebras of density character $\leq\aleph_\gamma$. 
\end{prop}

\begin{proof}
Let $K$ be the set of cardinals less than or equal to $\aleph_\gamma$. Then
$|K| = |\gamma|+\aleph_0$. By Proposition~\ref{P.UHF-classification}, the
number of isomorphism classes of UHF algebras of density character
$\leq\aleph_\gamma$ is equal to $|\{f : f\colon \cP\to K\}|=|K|^{\aleph_0}$.
\end{proof}

Note that for any ordinal $\gamma$ with $0\leq |\gamma|\leq
2^{\aleph_0}$, we have $(|\gamma|+\aleph_0)^{\aleph_0} = 2^{\aleph_0}$. Thus
for such $\gamma$, there are only as many UHF algebras of density character
$\leq\aleph_\gamma$ as there are separable UHF algebras (see the table in \S\ref{Sec:Open}).

\section{Non-classification of AM algebras in  regular uncountable character densities}
\label{Sec:Non-classAM-1}

The main result of this section shows that for a regular uncountable cardinal $\kappa$
there are as many AM algebras of
density character $\kappa$ as there are C*-algebras of density character
$\kappa$ and as many hyperfinite II$_1$ factors of density character $\kappa$ 
as there are II$_1$ factors of density character $\kappa$. 
The latter fact is in stark contrast with the separable case, when the hyperfinite~II$_1$ factor is unique. 
While there are continuum many separable UHF algebras, one should note that 
all AM algebras constructed here have the same K-theory as the (unique) CAR algebra. 

We first concentrate on case when $\kappa=\aleph_1$. 
Let $\Lambda$ be the set of all limit ordinals in $\aleph_1$. As an ordered
set, $\Lambda$ is isomorphic to $\aleph_1$. For each $\xi \in \aleph_1$, let
$A_\xi$ be the C*-algebra generated by two self-adjoint unitaries $v_\xi,
w_\xi$ with $v_\xi w_\xi = - w_\xi v_\xi$. By  \cite[Lemma~\ref{UHF1-Lem:M2}]{FaKa:Nonseparable}, $A_\xi$ is
isomorphic to $M_2(\bbC)$. We define a UHF algebra $A$ by $A := \bigotimes_{\xi \in
\aleph_1} A_\xi \cong \bigotimes_{\aleph_1} M_2(\bbC)$. For a subset $Y$ of
$\aleph_1$, we set $A_Y = \bigotimes_{\xi \in Y} A_\xi \subset A$. For $\xi \in
\aleph_1$, we use the notations $[0,\xi)$ and $[0,\xi]$ to denote the subsets
$\{\delta\in \aleph_1 : \delta<\xi\}$ and $\{\delta\in \aleph_1 : \delta \leq
\xi\}$ of $\aleph_1$. For each $\delta \in \Lambda$, we define $\alpha_\delta
\in \Aut(A)$ by
\[
\alpha_\delta = \bigotimes_{\xi \in [0,\delta)} \Ad v_\xi.
\]
Then we have $\alpha_\delta^2=\id$ and $\{\alpha_\delta\}_{\delta \in \Lambda}$
commute with each other. Let $G_\Lambda$ be the discrete abelian group of all
finite subsets of $\Lambda$ as in \cite[Definition~\ref{UHF1-Def:group}]{FaKa:Nonseparable}. Define an action
$\alpha$ of $G_\Lambda$ on $A$ by $\alpha_F := \prod_{\delta \in
F}\alpha_\delta$ for $F \in G_\Lambda$ and let $B := A \rtimes_{\alpha} G$. For
each $\delta \in \Lambda$, the unitary implementing $\alpha_\delta$ will be
denoted by $u_\delta \in B$. For a subset $S$ of $\Lambda$, we define $B_S :=
C^*(A \cup \{u_\delta\}_{\delta \in S}) \subset B$. We note that $B_S$ is
naturally isomorphic to $A \rtimes_{\alpha} G_S$ where $G_S$ is considered as a
subgroup of~$G_\Lambda$.

\begin{definition}
Let $S$ be a subset of $\Lambda$, and $\lambda$ be an element of $\Lambda$. We
define a subalgebra $D_{S,\lambda}$ of $B_{S}$ by
\[
D_{S,\lambda} := C^*\big(A_{[0,\lambda)} \cup \{u_{\delta}\}_{\delta \in S\cap
[0,\lambda)}\big)\subset B_{S}.
\]
\end{definition}

\begin{lemma}\label{Lem:B=limD}
For each $S \subset \Lambda$ the algebra $B_S$ is AM. 
Also,  $\{D_{S,\lambda}\}_{\lambda\in \Lambda}$ is a
$\sigma$-complete directed family  subalgebras of $B_S$ 
isomorphic to the CAR algebra
with dense union. 
\end{lemma}

\begin{proof} 
Consider a triple $(F,G,H)$ such that $F \subset \lambda$, $G
=\{\delta_1, \delta_2, \ldots, \delta_m\} \subset S$ and $H =\{\xi_1, \xi_2,
\ldots, \xi_m\} \subset \lambda$ are finite sets,  $F \cap H =
\emptyset$, and
\[
\xi_1 < \delta_1 < \xi_2 < \delta_2 < \xi_3 < \cdots < \delta_{m-1} < \xi_m <
\delta_m.
\]
For such  $(F,G,H)$ define $D_{(F,G,H)} \subset B_S$ by
\[
D_{(F,G,H)} := C^*\big(A_F \cup \{u_{\delta}\}_{\delta \in G} \cup
\{w_\xi\}_{\xi \in H}\big) \subset B_{S}.
\]
We have $A_F \cong M_{2^n}(\bbC)$ where $n$ is the cardinality of $F$. For each $k =
1, 2, \ldots , m$, there exists a unitary $v_k \in A_F$ with $v_k a =
\alpha_{\delta_k}(a) v_k$ for all $a \in A_F$. For $k = 1, 2, \ldots , m$, we
set $v_k' := v_ku_{\delta_k}$ which is a self-adjoint unitary in $D_{(F,G,H)}$
commuting  with $A_F$. We define self-adjoint unitaries $\{w'_k\}_{k=1}^m$ in
$D_{(F,G,H)}$ by $w'_k := w_{\xi_k}w_{\xi_{k+1}}$ for $k = 1, 2, \ldots , m -
1$ and $w'_m := w_{\xi_m}$. Since $F \cap H = \emptyset$, the unitaries
$\{w'_k\}_{k=1}^m$ commute with $A_F$. It is routine to check $v_k' w'_l = w'_l
v_k'$ for $k,l \in \{1, 2, \ldots , m\}$ with $k \neq l$, and $v_k' w'_k = -
w'_k v_k'$ for $k = 1, 2, \ldots , m$. Thus by  \cite[Lemma~\ref{UHF1-Lem:M2}]{FaKa:Nonseparable} the
subalgebra $A_k'$ of $D_{(F,G,H)}$ generated by $v_k'$ and $w'_k$ is isomorphic
to $M_2(\bbC)$ for every $k$. The family $\{A_F\} \cup \{A_k'\}_{k=1}^m$ of unital
subalgebras of $D_{(F,G,H)}$ mutually commutes, and generate $D_{(F,G,H)}$.
Hence $D_{(F,G,H)}$ is isomorphic to $M_{2^{n+m}}(\bbC)$.

For two such triples $(F,G,H),(F',G',H')$, we have $D_{(F,G,H)} \subsetneq
D_{(F',G',H')}$ if $F\cup H \subset F'$ and $G \subset G'$. Since there exist
infinitely many elements of $\lambda$ between two elements of $S$, for
arbitrary finite subsets $F \subset \lambda$ and $G \subset S$ there exists a
finite subset $H \subset \lambda$ such that the triple $(F,G,H)$ satisfies the
conditions above. Therefore the family $\{D_{(F,G,H)}\}_{(F,G,H)}$ of full
matrix subalgebras of $D_{S,\lambda}$ is directed. It is clear that the union
of this family is dense in $D_{S,\lambda}$. Since $D_{S,\lambda}$ is separable
and a unital direct limit of algebras $M_{2^k}(\bbC)$, $k\in \bbN$, it is
isomorphic to the CAR algebra.

Since the family $\{D_{S,\lambda}\}_{\lambda\in \Lambda}$ is clearly 
$\sigma$-complete and covers $B_S$,  
this completes the proof. 
\end{proof}

\begin{prop}\label{Prop:BSAM-1}
For every $S \subset \Lambda$, $B_S$ is a unital AM algebra of 
density character $\aleph_1$ with the same $K_0$, $K_1$, and the Cuntz semigroup as the
CAR algebra.
\end{prop}

\begin{proof}
Since $\chi(A)=\aleph_1$ and $|G_\Lambda|=\aleph_1$, $\chi(B_S)=\aleph_1$. By
Lemma~\ref{Lem:B=limD} the algebra $B_S$ is the direct limit of the
$\sigma$-complete system $D_{S,\lambda}$, $\lambda\in \Lambda$, of its
separable subalgebras each of which is isomorphic to the CAR algebra. By
Lemma~\ref{L.K-0.reflection} and \cite[Lemma~\ref{UHF1-Lem:phi=id}]{FaKa:Nonseparable}, $B_S$ has the same
$K_0$, $K_1$, and the Cuntz semigroup as the CAR algebra.
\end{proof}

\begin{lemma}\label{Lem:D'D''}
For $S \subset \Lambda$ and $\lambda \in \Lambda$,
we have
\begin{align*}
Z_{B_S}(D_{S,\lambda}) &= C^*\big(A_{\aleph_1\setminus [0,\lambda)} \cup
\{u_{\delta}u_{\delta'}\}_{\delta, \delta' \in S \setminus [0,\lambda)}
\big), \\
Z_{B_S}\big(Z_{B_S}(D_{S,\lambda})\big)
&= C^*\big(A_{[0,\lambda)}
\cup \{u_{\delta}\}_{\delta \in S\cap [0,\lambda]}\big).
\end{align*}
In particular,  $D_{S,\lambda} = Z_{B_S}(Z_{B_S}(D_{S,\lambda}))$ if and
only if $\lambda \notin S$.
\end{lemma}

\begin{proof}
Let us set $D':=C^*\big(A_{\aleph_1\setminus [0,\lambda)} \cup
\{u_{\delta}u_{\delta'}\}_{\delta, \delta' \in S \setminus [0,\lambda)} \big)$.
It is clear that $A_{\aleph_1\setminus [0,\lambda)} \subset
Z_{B_S}(D_{S,\lambda})$ and $u_{g} \in Z_{B_S}(D_{S,\lambda})$ for $g \in G_S$
such that $|g|$ is even and $g \subset [\lambda,\aleph_1)$. Hence we get $D'
\subset Z_{B_S}(D_{S,\lambda})$. Take $a \in Z_{B_S}(D_{S,\lambda})$. For any
$\e > 0$, there exist a finite set $F\subseteq \aleph_1$, finite families
$b_1,b_2,\ldots, b_n \in A_F$ and $g_1, g_2, \ldots, g_n \in G_S$ such that $b
= \sum_{k=1}^n b_k u_{g_k}$ satisfies $\|a - b \| < \e$. Let $\delta_1,
\delta_2, \ldots, \delta_m$ be the list of $[0,\lambda)\cap (\bigcup_{k=1}^n
g_k)$ ordered increasingly. Choose $H = \{\xi_1, \xi_2, \ldots, \xi_m\} \subset
\aleph_1 \setminus F$ such that
\[
\xi_1 < \delta_1 < \xi_2 < \delta_2 < \xi_3 < \cdots
< \delta_{m-1} < \xi_m < \delta_m < \lambda.
\]
For each $H'\subset H$, we define a self-adjoint unitary $w_{H'}$ by $w_{H'} =
\prod_{\xi \in H'}w_{\xi}$ Let us define a linear map $E \colon B_S \to B_S$ by
$E(x) = 2^{-m}\sum_{H'\subset H}w_{H'}xw_{H'}$. Then $E$ is a contraction.
Since $a \in Z_{B_S}(D_{S,\lambda})$, we have $E(a)=a$. Hence $\|a - E(b)\| <
\e$. For $g \in G_S$ with $\delta_k \in g$ for some $k$, we have $E(u_g)=0$.
For $g \in G_S$ such that $g \subset [\lambda,\aleph_1)$ and $|g|$ is odd, we
also have $E(u_g)=0$. For $g \in G_S$ with $g \subset [\lambda,\aleph_1)$ and
$|g|$ is even, we get $E(u_g)=u_g$. Therefore $E(b) = \sum_{k} b_k u_{g_k}$
where $k$ runs over elements such that $g_k \in G_S$ satisfies that $|g_k|$ is
even and $g_k \subset [\lambda,\aleph_1)$. Next let $F' = F \cap [0,\lambda)$.
We define a contractive linear map $E' \colon B_S \to B_S$ by $E'(x) =
\int_{U}uxu^* du$ where $U$ is the unitary group of the finite dimensional
subalgebra $A_{F'}$ of $D_{S,\lambda}$, and $du$ is its normalized Haar
measure. Since $a \in Z_{B_S}(D_{S,\lambda})$, we have $E'(a)=a$. Hence $\|a -
E'(E(b))\| < \e$. For $g_k \in G_S$ such that $|g_k|$ is even and $g_k \subset
[\lambda,\aleph_1)$, we have $u_{g_k}u^* =u^*u_{g_k}$ for all $u \in U$. Hence
for such $k$, we have $E'(b_k u_{g_k})=E'(b_k) u_{g_k}$. Since $E'(b_k) \in
A_{[\lambda,\aleph_1)}$, we get $E'(E(b)) \in D'$. Since $\e >0 $ was
arbitrary, $a \in D'$. Thus we have shown $Z_{B_S}(D_{S,\lambda}) = D'$.

The equality $Z_{B_S}(D') = C^*\big(A_{[0,\lambda)} \cup \{u_{g}\}_{\text{$g
\in G_S$, $g \subset [0,\lambda]$}}\big)$ can also be proved in a similar way
as above. The only difference is that $\delta_1, \delta_2, \ldots, \delta_m$ is
now the list of $(\lambda,\aleph_1)\cap (\bigcup_{k=1}^n g_k)$ ordered
increasingly, and choose $H = \{\xi_1, \xi_2, \ldots, \xi_m\} \subset \aleph_1
\setminus F$ such that
\[
\lambda < \xi_1 < \delta_1 < \xi_2 < \delta_2 < \xi_3 < \cdots
< \delta_{m-1} < \xi_m < \delta_m.
\]
We leave the details to the readers.
\end{proof}

\begin{lemma}\label{Lem:B=C*(D,D')}
For $S \subset \Lambda$ and $\lambda \in \Lambda$,
$B_S$ is generated by $D_{S,\lambda}$ and
$Z_{B_S}(D_{S,\lambda})$
if and only if
$S \subset [0,\lambda)$.
\end{lemma}

\begin{proof}
Lemma~\ref{Lem:D'D''} implies that $B_S$ is generated by $D_{S,\lambda}$ and
$Z_{B_S}(D_{S,\lambda})$ if $S \subset [0,\lambda)$.
If there exists $\delta \in S \setminus [0,\lambda)$,
then $u_\delta$ is not in the C*-algebra
generated by $D_{S,\lambda}$ and $Z_{B_S}(D_{S,\lambda})$.
\end{proof}

Compare the following proposition to Proposition~\ref{P.forcing}. 

\begin{prop}
For $S \subset \Lambda$,
the C*-algebra $B_S$ is UHF if and only if $S$ is bounded.
In this case, $B_S$ is isomorphic to $A \cong \bigotimes_{\aleph_1} M_2(\bbC)$.
\end{prop}

\begin{proof}
When $S$ is unbounded, the $\sigma$-complete system
$\{D_{S,\lambda}\}_{\lambda\in \Lambda}$ in Lemma~\ref{Lem:B=limD} satisfies
that $B_S$ is not generated by $D_{S,\lambda}$ and $Z_{B_S}(D_{S,\lambda})$ for
all $\lambda$ by Lemma~\ref{Lem:B=C*(D,D')}. Hence $B_S$ is not a UHF algebra.
When $S \subset [0, \lambda)$ for some $\lambda\in\aleph_1$, then we have $B_S
= D_{S,\lambda} \otimes A_{[\lambda,\aleph_1)}$ by Lemma~\ref{Lem:B=C*(D,D')}.
By Lemma~\ref{Lem:B=limD}, $D_{S,\lambda}$ is the CAR algebra. Hence $B_S \cong
\bigotimes_{\aleph_1} M_2(\bbC)$.
\end{proof}

\begin{prop}\label{P.S-S'}
Let $S$ and $S'$ be two subsets of $\Lambda$.
If $B_S$ and $B_{S'}$ are isomorphic,
then there exists a club $\Lambda_0$ in $\Lambda$
such that $\Lambda_0\cap (S\Delta S')=\emptyset$.
\end{prop}

\begin{proof}
Assume $\Phi\colon B_S\to B_{S'}$ is an isomorphism. By
\cite[Proposition~\ref{UHF1-Prop:find_club}]{FaKa:Nonseparable}, there exists a club $\Lambda_0 \subset
\Lambda$ such that $\Phi[D_{S,\lambda}]=D_{S',\lambda}$ for all $\lambda \in
\Lambda_0$. For $\lambda \in \Lambda_0$, $\lambda \in S$ if and only if
$\lambda \in S'$ by Lemma~\ref{Lem:D'D''}. Thus we have $\Lambda_0\cap (S\Delta
S')=\emptyset$.
\end{proof}

\begin{proof}[Proof of Theorem~\ref{T.non-classification}]
By Lemma~\ref{L.many-subsets} we can fix a family $S_0(X)$, $X\subseteq
\aleph_1$, of subsets of $\aleph_1$  such that $S_0(X)\Delta S_0(Y)$ is
stationary whenever $X\neq Y$. Since $\Lambda$ is a club in $\aleph_1$, the
sets $S(X)=\Lambda\cap S_0(X)$ retain this property.

Therefore the algebras  $B_{S(X)}$, $X\subseteq \aleph_1$, are nonisomorphic by
Proposition~\ref{P.S-S'}. By Proposition~\ref{Prop:BSAM-1} these algebras have
the same $K$-theory and Cuntz semigroup as the CAR algebra.
\end{proof}

For any uncountable regular cardinal $\kappa$ one can define $<\kappa$-complete directed
systems of algebras of density character $<\kappa$ and prove results analogous
to those for $\sigma$-complete directed systems so that the latter coincide
with   $<\aleph_1$-complete systems. Given this and Lemma~\ref{L.many-subsets},
a straightforward extension of the proof of Theorem~\ref{T.non-classification}
gives the following.

\begin{theorem} \label{T.non-classification.2}
If $\kappa$ is a regular cardinal then there are $2^\kappa$
nonisomorphic AM algebras of density character $\kappa$. \qed
\end{theorem}

However, this method does not work for singular cardinals and we shall treat this case in the following section.

\section{Non-classification of AM algebras in all character densities}
\label{Sec:Non-classAM-2}

The proof of the present section relies on two components. The first is the  
 non-structure theory as developed by Shelah 
in \cite{Sh:E59} and adapted to metric structures in \cite{FaSh:Dichotomy}, 
and the second is the order property of theories of C*-algebras and II$_1$ factors
proved in \cite{FaHaSh:Model1}.   
Readers with  background in model theory will notice that 
the algebras that we construct are EM-models 
generated by indiscernibles which 
witness that their theory has the order property.

Fix a total ordering  $\Lambda$ and let $\Lambda^+$ denote $\Lambda\times \bbN$
with the lexicogaphical ordering. We identify $\Lambda$ with $\Lambda\times\{0\}\subseteq \Lambda^+$ and note that between any two elements $\xi<\eta$  of $\Lambda$ there are infinitely many elements of $\Lambda^+\setminus \Lambda$. 
For each $\xi \in \Lambda^+$, let
$A_\xi$ be the C*-algebra generated by two self-adjoint unitaries $v_\xi,
w_\xi$ with $v_\xi w_\xi = - w_\xi v_\xi$. By \cite[Lemma~\ref{UHF1-Lem:M2}]{FaKa:Nonseparable}, $A_\xi$ is
isomorphic to $M_2(\bbC)$. We define a UHF algebra $A_\Lambda$ 
by $A_\Lambda := \bigotimes_{\xi \in
\Lambda^+} A_\xi \cong \bigotimes_{\Lambda^+} M_2(\bbC)$.  For $\xi \in
\Lambda^+$ we write 
\begin{align*}
[0,\xi):=& 
\{\delta\in \Lambda^+ : \delta<\xi\}\\
 [0,\xi]:=&\{\delta\in \Lambda^+: \delta \leq
\xi\}. 
\end{align*} 
 For each $\delta \in \Lambda$, we define $\alpha_\delta
\in \Aut(A)$ by
\[
\alpha_\delta =\textstyle \bigotimes_{\xi \in [0,\delta)} \Ad v_\xi.
\]
Then we have $\alpha_\delta^2=\id$ and $\{\alpha_\delta\}_{\delta \in \Lambda}$
commute with each other. Let $G_\Lambda$ be the discrete abelian group of all
finite subsets of $\Lambda$ as in 
\cite[Definition~\ref{UHF1-Def:group}]{FaKa:Nonseparable}. 
Define an action
$\alpha$ of $G_\Lambda$ on $A_\Lambda$ by $\alpha_F := \prod_{\delta \in
F}\alpha_\delta$ for $F \in G_\Lambda$ and let $B_\Lambda := A_\Lambda \rtimes_{\alpha} G_\Lambda$. For
each $\delta \in \Lambda$, the unitary implementing $\alpha_\delta$ will be
denoted by $u_\delta \in B$. 
For $S\subseteq \Lambda$ let $A_S:=\bigotimes_{\xi\in S\times \bbN} A_\xi$
and consider it as a subalgebra of $A_\Lambda$.

\begin{definition}
If $S$ is a subset of $\Lambda$
define a subalgebra $D_{S}$ of $B_\Lambda$ by
\[
D_{S} := C^*\big(A_{S} \cup \{u_{\delta}\}_{\delta \in S}\big).
\]
\end{definition}

\begin{lemma}\label{Lem:B=limD-2}
For each uncountable total order $\Lambda$ the algebra $A_\Lambda$ is AM.  
Also, 
 $\{D_{S}:S\in [\Lambda]^{\aleph_0}\}$ is a
$\sigma$-complete directed family  of subalgebras of $A_\Lambda$ 
isomorphic to the CAR algebra
with dense union. 
\end{lemma}

\begin{proof} This proof is almost identical to the proof of Lemma~\ref{Lem:B=limD}. 
The assumption that $\lambda$ is a limit ordinal used in the former proof is 
replaced by the fact that the generators are indexed by $\Lambda^+$. 
\end{proof} 

\begin{prop}\label{Prop:BSAM}
For every infinite cardinal $\kappa$ and 
total ordering $ \Lambda$ of cardinality $\kappa$, 
$B_\Lambda$ is a unital AM algebra of character
density equal to $\kappa$ 
 with the same $K_0$, $K_1$, and the Cuntz semigroup as the
CAR algebra.
\end{prop}

\begin{proof} 
Since $\chi(A_\Lambda)=\kappa$ and $|G_\Lambda|=\kappa$, $\chi(B_\Lambda)=\kappa$. By
Lemma~\ref{Lem:B=limD-2} the algebra $B_S$ is the direct limit of the
$\sigma$-complete system $D_{S}$, $S\in [\Lambda]^{\aleph_0}$, of its
separable subalgebras each of which is isomorphic to the CAR algebra. By
Lemma~\ref{L.K-0.reflection} and \cite[Lemma~\ref{UHF1-Lem:phi=id}]{FaKa:Nonseparable}, 
$B_S$ has the same
$K_0$, $K_1$, and the Cuntz semigroup as the CAR algebra.
\end{proof}

Assume $P(\vec x, \vec y)$ is a *-polynomial in $2n$ variables. 
Then for every C*-algebra $A$ the expression 
$\phi(\vec x, \vec y)=\|P(\vec x, \vec y)\|$
defines a uniformly continuous map from $A^{2n}$ into the nonnegative reals. 
Let $(A_{\leq 1})$ denote the unit ball of~$A$ and on $(A_{\leq 1})^n$ define
a binary relation $\prec_\phi$ by 
letting $\vec a\prec \vec b$ if 
\[
\phi(\vec a, \vec b)=1 \qquad\text{ and } 
\qquad \phi(\vec b, \vec a)=0. 
\]
Note that $\prec_\phi$ is not required to be an ordering. 
If $\Lambda$ is a total ordering   we shall say that
an indexed set $\vec a_\lambda$, for $\lambda\in \Lambda$
is a \emph{$\phi$-chain} 
if $\vec a_\lambda\prec_\phi \vec a_{\lambda'}$ whenever $\lambda<\lambda'$. 
We write $\vec a\preceq_\phi \vec b$ if $\vec a=\vec b$ or $\vec a\prec_\phi \vec b$.

\begin{definition}[{\cite[Definition~3.1]{FaSh:Dichotomy}}] \label{Def.wsl}
A $\phi$-chain $\cC$ is \emph{weakly $(\aleph_1,\phi)$-skeleton like} inside $A$ if for every 
$\vec a\in A^n$ there is a countable $\cC_{\vec a}\subseteq \cC$ such 
that for all $\vec b$ and $ \vec c$ in $\cC$ for which  we have $\vec b\preceq_\phi \vec c$ and  
no $\vec d\in \cC_{\vec a}$ satisfies $\vec b\preceq_\phi \vec d\preceq_\phi \vec c$ we have
\[
 \phi(\vec b,\vec a)=\phi(\vec c,\vec a)
\qquad 
 \text{and} 
 \qquad
 \phi(\vec a,\vec b)=\phi (\vec a,\vec c).
\]
\end{definition}

\begin{lemma} \label{L.nonstructure}
Assume $\cK$ is a class of C*-algebras,  $\phi(\vec x, \vec y)$ is as above, 
and~$\kappa$ is an uncountable cardinal. 
If for every linear ordering $\Lambda$ of cardinality $\kappa$ 
there is $B_\Lambda\in \cC$ of density character $\kappa$ 
such that  the $n$-th power of the unit ball of $B_\Lambda$  includes
 a $\phi$-chain $\cC$ isomorphic to $\Lambda$ which  is weakly $(\aleph_1,\phi)$-skeleton like, 
 then $\cK$ contains $2^\kappa$ nonisomorphic algebras of density character~$\kappa$.  
\end{lemma}

\begin{proof} This is an immediate consequence of  results from \cite{FaSh:Dichotomy}, but we sketch a proof for the convenience of the reader. 
By \cite[Lemma 2.5]{FaSh:Dichotomy} for every $m\in \bbN$ (actually $m=3$ suffices)
there are~$2^\kappa$ total orderings
of cardinality $\kappa$ that have disjoint representing sequences of $m,\kappa$-invariants
(in the sense of \cite[\S 2.2]{FaSh:Dichotomy}). 
 For any such ordering~$\Lambda$ the algebra $B_\Lambda$ has 
 density character $\kappa$ and therefore the  $m,\kappa$-invariant of $\Lambda$ belongs to  
  $\INV^{m,\kappa}(B_\Lambda)$, as defined in \cite[Definition 3.8 and \S 6.2]{FaSh:Dichotomy}.
By \cite[Lemma~6.4]{FaSh:Dichotomy} for each C*-algebra $B$ of density character 
$\kappa$ the set $\INV^{m,\kappa}(B)$ has cardinality at most~$\kappa$. 
Since $2^\kappa$ cannot be written as the supremum 
of $\kappa$ smaller cardinals (\cite[Corollary~10.41]{Ku:Book}), 
 by a counting argument there are $2^\kappa$ isomorphism classes among algebras~$B_\Lambda$ for a total ordering  $\Lambda$ of cardinality $\kappa$.  
\end{proof} 

\begin{proof}[Proof of Theorem~\ref{T.non-classification}]
Formula
$
\phi(x_1,x_2,y_1,y_2)=\frac 12
\|[x_1,y_2]\|
$
defines a uniformly continuous function on $A^4$ for any C*-algebra $A$. 
With $\Lambda$, 
 $B_\Lambda$, $u_\xi$, and $w_\xi$ as in the first paragraph of \S\ref{Sec:Non-classAM-1}, 
   for all $\xi$ and $\eta$ in $\Lambda$ we have
\[
\phi(u_\xi, w_\xi, u_\eta, w_\eta)=
\begin{cases} 0, \text{ if }\xi<\eta\\
1, \text{ if } \xi\geq \eta, 
\end{cases}
\]
and therefore $(u_\xi,w_\xi)$, for $\xi\in \Lambda$, is a  $\phi$-chain. 

Consider $S\subseteq \Lambda$. 
 On the set $\Lambda\setminus S$ define 
an equivalence relation, $\xi\sim_S \eta$ if and only if no element of $S$ is 
between $\xi$ and $\eta$. Then for $\xi\sim_S \eta$ we have that 
the algebras $C^*(B_S\cup\{u_\xi,w_\xi\})$ and $C^*(B_S\cup \{u_\eta,w_\eta\})$ 
are isomorphic via an isomorphism that is an identity on $B_S$ and 
sends $u_\xi$ to $u_\eta$ and $w_\xi$~to~$w_\eta$.

We claim  that   
$\Lambda$ is weakly $(\aleph_1,\phi)$-skeleton like in $B_S$. 
First note that  every finite set  $F\subseteq  B_\Lambda$ is included in $D_S$ for some 
 countable $S=S(F)\subseteq \Lambda$.   
 For $a_1$ and $a_2$ in $B_\Lambda$ 
 fix  a countable $S$  such that $\{a_1,a_2\}\subseteq B_S$. 
 Then let ${\mathcal C}_{\{a_1,a_2\}}=S$
 and note that $\xi\sim_S \eta$ implies that  
$\phi(a_1,a_2,u_\xi,w_\xi)=\phi(a_1,a_2,u_\eta,w_\eta)$
and $\phi(u_\xi,w_\xi,a_1,a_2)=\phi(u_\eta,w_\eta,a_1,a_2)$.

Therefore our 
 distinguished $\Lambda$-chain $(u_\xi,w_\xi)$, for $\xi\in \Lambda$, 
 s $(\aleph_1,\phi)$-skeleton like  
 Lemma~\ref{L.nonstructure} applies to show that there are $2^\kappa$ isomorphism classes
among algebras $B_\Lambda$ for $|\Lambda|=\kappa$. 
 By Proposition~\ref{Prop:BSAM} these algebras have
the same $K$-theory and Cuntz semigroup as the CAR algebra.
\end{proof}

The assumption that we were dealing with C*-algebras in  Lemma~\ref{L.nonstructure} 
was not crucial. This lemma 
 applies to any class of models of logic of metric structures (\cite{BYBHU}, \cite{FaHaSh:Model2}), 
 and in particular to II$_1$ factors. 
We shall now state the general form of  Lemma~\ref{L.nonstructure}. 
The definition of
`metric structure' and `formula' is given  in \cite{BYBHU} (see also \cite{FaHaSh:Model2} for 
the case of C*-algebras and tracial von Neumann algebras). Although this lemma 
uses logic for metric structures, we note that class $\cC$ is not required to be axiomatizable. 
Indeed, neither AM algebras nor hyperfinite II$_1$ factors are axiomatizable 
(cf. the proof of  \cite[Proposition~6.1]{FaHaSh:Model2}, but see also 
\cite{Mitacs2012}). We state this lemma in the case
of bounded metric structures, and the version for II$_1$ factors necessitates requiring that 
the chain  be included in the $n$-power of the unit ball. The proof of Lemma~\ref{L.nonstructure+} 
 is identical to the proof of 
Lemma~\ref{L.nonstructure}. 

\begin{lemma} \label{L.nonstructure+}
Assume $\cC$ is a class of bounded metric structures and $\phi(\vec x, \vec y)$ is a $2n$-ary formula. 
Assume that 
 for every linear ordering $\Lambda$ there is $A_\Lambda\in \cC$ of density character $|\Lambda|$ 
such that  $A_\Lambda^n$  includes
 a $\phi$-chain $\cC$ isomorphic to~$\Lambda$ which  is weakly $(\aleph_1,\phi)$-skeleton like.
 Then $\cC$ contains $2^\kappa$ nonisomorphic structures in every uncountable 
 density character $\kappa$. \qed
\end{lemma}

\begin{proof}[Proof of Theorem~\ref{T.non-classification.hyperfinite}]
For each of the algebras $B_\Lambda$ constructed in the proof of Theorem~\ref{T.non-classification} consider the GNS representation corresponding to its unique trace
and let $R_\Lambda$ be the weak closure of the image of $B_\Lambda$. 
Then each $R_\Lambda$ is a hyperfinite II$_1$ factor whose predual has density character 
$\kappa=|\Lambda|$. 
The formula
\[
\psi(x_1,x_2,y_1,y_2)=\|[x_1,y_2]\|_2
\]
defines a uniformly continuous with respect to the 2-norm function on $R_\Lambda$. 
Let  $A_\Lambda$ denote the operator norm unit ball of $B_\Lambda$. Then each $A_\Lambda$, 
equipped with the $\ell_2$-norm and function that evaluates $\psi$  is a
 bounded metric structure  and it suffices to check that 
 Lemma~\ref{L.nonstructure+} applies to this family. 
Again $(u_\xi,w_\xi)$, for $\xi\in \Lambda$, is  
a $\phi$-chain that is weakly $(\aleph_1,\phi)$-skeleton like.
By  Lemma~\ref{L.nonstructure+} 
there are $2^\kappa$ nonisomorphic unit balls of II$_1$-factors of the form $B_\Lambda$
with the predual of density character~$\kappa$. Therefore there are $2^\kappa$ nonisomorphic
hyperfinite II$_1$ factors with the density character $\kappa$ for every uncountable cardinal $\kappa$.   
\end{proof}

Note that the assumption that $\kappa$ is uncountable is necessary in Lemma~\ref{L.nonstructure}, since 
the hyperfinite II$_1$ factor with separable predual is unique.

The remainder of this section is aimed at logicians. 
A class of models is non-classifiable in a strong sense if it does not allow 
 sequences of cardinal numbers as complete invariants (see \cite{Sh:Classification}). 
 The following proposition shows that AM algebras and hyperfinite II$_1$ factors 
are non-classifiable even in this strong sense.

\begin{prop} \label{P.forcing}
 There are AM algebras $A$ and $B$ of density character~$\aleph_1$ 
 and a forcing notion $\bbP$ that does not collapse cardinals or add countable sequences of cardinals  such that $A$ and $B$ are not isomorphic, but $\bbP$ forces that $A$ and $B$ 
 are isomorphic. 
 
 There are also hyperfinite II$_1$ factors of density character $\aleph_1$ with the same property. 
 \end{prop} 
 
 \begin{proof} Let $S\subseteq \aleph_1$ a stationary set whose complement is also stationary. 
 Let $\eta$ denote the ordering of the rational numbers. 
 Let $\Lambda(1)$ be the linear ordering obtained from $\aleph_1$ by replacing
 all points  with a copy of $\eta$ (i.e., the lexicographical ordering of $\aleph_1\times \eta$). 
  Let $\Lambda(2)$ be the linear ordering obtained from $\aleph_1$ by replacing
 points in $S$ with a copy of $\eta$ and leaving points in $\aleph_1\setminus S$ unchanged. 
 
Since $\aleph_1\setminus S$ is stationary,  
the argument from the proof of  Proposition~\ref{P.S-S'}
shows that the   algebras $A:=A_{\Lambda(1)}$ and $B:=A_{\Lambda(2)}$ 
are not isomorphic. 
Let $\bbP$ be Jensen's forcing for adding  a club subset of  $S$. 
Then $\bbP$ is $\sigma$-distributive  (see e.g., \cite[VII.H18]{Ku:Book})
and therefore it does not collapse $\aleph_1$ and does not add 
 new sequences of cardinals. 
Since $\bbP$ has cardinality $\aleph_1$, it does not collapse cardinals larger than $\aleph_1$. 

We claim that $\bbP$ nevertheless forces $A$ and $B$ to be isomorphic. 
It clearly suffices to show that it forces $\Lambda(1)$ and $\Lambda(2)$ are isomorphic
linear orderings. 
If $C\subseteq S$ is the club added by $\bbP$, then points in $C$ separate $\Lambda(1)$ and $\Lambda(2)$ into $\aleph_1$ sequence of countable linear orderings without endpoints. 
Any two such orderings are isomorphic by Cantor's classical back-and-forth argument, 
and these isomorphisms  together define an isomorphism between $\Lambda(1)$ and $\Lambda(2)$. 

Construction of the required II$_1$ factors is analogous. 
\end{proof} 

 Proposition~\ref{P.forcing} shows that the   classification
problem of AM algebras of density character $\aleph_1$ is at least as
complicated as the classification of subsets of $\aleph_1$ modulo the
nonstationary ideal. The latter problem is largely considered to be
intractable.

\section{Concluding remarks}\label{Sec:Open}

The number of AM algebras and 
 UHF algebras in some character densities 
as well as the number of  hyperfinite II$_1$-factors whose predual has the same density character  
is given in the  table below.  We  identify each 
cardinal with the least ordinal having it as a cardinality, 
write  $\fc:=2^{\aleph_0}$, and $\fc^+$ denotes the least cardinal 
greater than $\fc$. 

\begin{table}[h]
 {\small\tt \begin{tabular}{r|cccccccccc}
\parbox{1.1in}{density\\ character\\} & $\aleph_0$ & $\aleph_1$ & $\aleph_2$& \dots 
& $\aleph_\omega$ & \dots & $\aleph_{\omega_1}$ &\dots & $\aleph_{\fc^+}$& \dots\\ 
\hline
\parbox{1.1in}{the number of\\ UHF  algebras\\ } & $\fc$ & $\fc$ & $\fc$ & \dots & $\fc$ & \dots & 
$\fc$ & \dots & $\fc^+$ & \dots  \\
\parbox{1.1in}{the number of\\ AM algebras\\ } & $\fc$ & $2^{\aleph_1}$ & $2^{\aleph_2}$ & \dots & $2^{\aleph_{\omega}}$ 
& \dots & 
$2^{\aleph_{\omega_1}}$ & \dots & $2^{\aleph_{\fc^+}}$ & \dots\\
\parbox{1.1in}{the number of\\ 
hyperfinite\\ II$_1$ factors\\ }   & $\fc$ & $2^{\aleph_1}$ & $2^{\aleph_2}$ & \dots & $2^{\aleph_{\omega}}$ 
& \dots & 
$2^{\aleph_{\omega_1}}$ & \dots & $2^{\aleph_{\fc^+}}$ & \dots
\end{tabular}}
\end{table} 

  While the cardinals in this table resemble those predicted by Shelah's 
Main Gap Theorem for the number of models of 
classifiable and non-classifiable theories in uncountable cardinalities (\cite{Sh:Classification}), 
it should be noted that all algebras appearing in our
proofs are elementarily equivalent to the CAR algebra and that the class of UHF algebras 
does not seem to have a natural model-theoretic characterization. 
On the other hand,  all AM (and even all LM) algebras are by \cite{Mitacs2012}, atomic models. 
It is not difficult to see that the methods of \cite{Mitacs2012} also show that  
hyperfinite II$_1$ factors (of arbitrary density character) are atomic models.

For a unital C*-algebra $A$, we define two generalized integers, $\sn(A)$ and
$\sntoo(A)$, associated to $A$ as follows. Recall that $\cP$ denotes the set of
all primes. If $A=\bigotimes_{p\in \cP}\bigotimes_{\kappa_p} M_p(\bbC)$, then
let $\sn(A)_p=\kappa_p$ for $p\in \cP$.
\begin{align*}
\sntoo(A)_p := \sup\big\{|X| : &\text{there exists a unital homomorphism}\\
&\text{from $\textstyle\bigotimes_X M_p(\bbC)$ to $A$}\big\}
\end{align*}
 for each $p \in
\mathcal{P}$. Clearly these two definitions coincide when $A$ is separable and
$\sn(A)\leq \sntoo(A)$, pointwise.

\begin{problem} \label{Prob.sn} If $A$ is a UHF algebra, is $\sn(A)=\sntoo(A)$?
\end{problem}

Here is a version of Problem~\ref{Prob.sn}.

\begin{problem} \label{Prob.emb} Assume $A$ is UHF,
$\kappa<\kappa'$ are
cardinals
 and $\bigotimes_{\kappa'}M_2(\bbC)$ unitally embeds into
$\bigotimes_{\kappa}M_2(\bbC)\otimes A$. Can we conclude that there is a unital
embedding of $\bigotimes_\kappa M_2(\bbC)$ into $A$?
\end{problem}

We cannot even prove that in the above situation $M_2(\bbC)$ unitally embeds
into~$A$. The most embarrassing version of Problem~\ref{Prob.emb} is whether
$\bigotimes_{\aleph_1}M_2(\bbC)$ unitally embeds into $\bigotimes_{\aleph_0}
M_2(\bbC)\otimes\bigotimes_{\aleph_1}M_3(\bbC)$.   Since any two unital copies
of $M_n(\bbC)$, for $n\in\bbN$, in a UHF algebra are conjugate,
Problem~\ref{Prob.emb} has a positive answer when $\kappa$ is finite. 

Standard results on classification of unital, separable, nuclear, simple C*-algebras  imply 
that if $A$ is not UHF then the answer to Problem~\ref{Prob.emb} is
negative.
We shall need the following well-known fact.

\begin{lemma}\label{L.group} There are an  abelian group $G$ and  a nonzero $g_0\in
G$ such that
\begin{enumerate}
\item \label{L.group.1} no infinite nonzero sequence $(f_n)$ in $G$ satisfies
$f_n=2f_{n+1}$ for all~$n$,
\item \label{L.group.2} for every $n$ there is $h_n$ satisfying $2^n h_n=g_0$, and
\item \label{L.group.3} $h_1\neq 2f$ for all $f\in G$.
\end{enumerate}
\end{lemma}

\begin{proof}
Let $m$ be the generalized integer defined by $m_2={\aleph_0}$, $m_p=0$ for
$p\neq 2$ and let $G$ be the subgroup of (recall that $\bbZ[1/m]$ was defined before Proposition~\ref{P.tracial-state}) $\prod_{n=1}^\infty
(\bbZ/2^n\bbZ)\times \bbZ[1/m]$ consisting of all $(x_n)_{n\leq\infty}$
such that  
\[
x_n=2^n x_\infty\mod 1
\]
 for all but finitely many $n$. Since $x_\infty\in \bbZ[1/m]$, for a large enough $n$ we will have 
 that $x_\infty$ is equal to an element of $\bbZ/2^n\bbZ$ for all large enough $n$. 
 Let
$g_0=(x_n)$ where $x_n=0$ for all finite $n$ and $x_\infty=1$. Then $h_n=(x_j)$
where $x_j=0$ for $j<n$, $x_j=2^{j-n}$ for $j\geq n$ finite and
$x_\infty=2^{-n}$ satisfy $2^nh_n=g_0$ for each~$n$.
If  $f=(y_n)$ is such that $2f=h_1$ then necessarily $2y_1=1$,  but there is no
such element in $\bbZ/2\bbZ$. This proves \eqref{L.group.3} and the proof of
\eqref{L.group.1} is similar.
\end{proof}

The following two propositions rely on the Kirchberg--Phillips classification
of Kirchberg algebras $A$ by its $K$-theoretic invariants
$$
I(A)=(K_0(A), [1]_0,
K_1(A))
$$
and the fact that every pair of countable abelian groups with a
distinguished element is an invariant of some Kirchberg algebra (see e.g.,
\cite{Phi:Classification} or \cite{Ror:Classification}).

\begin{prop}\label{L.sn(A)}
There is a C*-algebra $A$ such that $\bigotimes_X M_2(\bbC)$ unitally embeds
into $A$ if and only if $X$ is finite. Moreover, $A=M_2(\bbC)\otimes B$ for
some $B$ such that $M_2(\bbC)$ does not unitally embed into $B$.
\end{prop}

\begin{proof}
 Let $G$ and $g_0$ be as in Lemma~\ref{L.group} and let $A$ be the Kirchberg algebra with $K_0(A)=G$ such
 that $g_0$ is equal to the class of the identity and with trivial $K_1(A)$.
For $n\in \bbN$, $M_{2^n}(\bbC)$ embeds unitally into $A$ by $K$-theoretic
consideration.
 Pick a projection $q$ in $A$ such that $[q]=g$, and let
$C$ be a unital copy of $M_2(\bbC)$ in $A$ with $q$ as its matrix unit. Then
$A\cong M_2(\bbC)\otimes B$, with $B=Z_A(C)$. By \eqref{L.group.3} and
$K$-theoretic considerations $B$ has no unital copy of~$M_2(\bbC)$ and the
proof is complete.
\end{proof}

\begin{prop}\label{P.5.7}
There is no unital *-homomorphism from $M_2(\bbC)$ into the Cuntz algebra
 ${\mathcal O}_3$, but there is a
unital *-homomorphism from the CAR algebra into $M_2({\mathcal O}_3)$.
\end{prop}

\begin{proof} Let $A$ denote the CAR algebra. The algebras $\cO_3$ and  $M_2(\cO_3)$
 are Kirchberg algebras. Since $I(\cO_3)=(\bbZ/2\bbZ,1,0)$ (see
\cite{Ror:Classification}), the identity in $K_0(\cO_3)$ is not divisible by 2
and therefore $M_2(\bbC)$ is not a unital subalgebra of $\cO_3$. Since
$M_2(\cO_3)\otimes\cK\cong \cO_3\otimes \cK$ we have  $K_0(M_2(\cO_3)=\bbZ/2\bbZ$ but the
class of the identity element is $0$ and we have
$I(M_2(\cO_3)=(\bbZ/2\bbZ,0,0)$.

Since $2\times 0=0$ and $M_2(\cO_3)$ is purely infinite, it 
contains a unital copy of $\cO_2$ and therefore a
unital copy of any other simple nuclear C*-algebra---including the CAR algebra. 
\end{proof}

\bibliographystyle{amsplain}
\bibliography{uhfbib}

\providecommand{\bysame}{\leavevmode\hbox to3em{\hrulefill}\thinspace}
\providecommand{\MR}{\relax\ifhmode\unskip\space\fi MR }
\providecommand{\MRhref}[2]{%
  \href{http://www.ams.org/mathscinet-getitem?mr=#1}{#2}
}
\providecommand{\href}[2]{#2}
\begin{thebibliography}{10}

\bibitem{BYBHU}
I.~Ben~Yaacov, A.~Berenstein, C.W. Henson, and A.~Usvyatsov, \emph{Model theory
  for metric structures}, Model Theory with Applications to Algebra and
  Analysis, Vol. II (Z.~Chatzidakis et~al., eds.), Lecture Notes series of the
  London Math. Society., no. 350, Cambridge University Press, 2008,
  pp.~315--427.

\bibitem{Black:Operator}
B.~Blackadar, \emph{Operator algebras}, Encyclopaedia of Mathematical Sciences,
  vol. 122, Springer-Verlag, Berlin, 2006, Theory of $C\sp *$-algebras and von
  Neumann algebras, Operator Algebras and Non-commutative Geometry, III.

\bibitem{BrOz}
N.~Brown and N.~Ozawa, \emph{{$C\sp *$}-algebras and finite-dimensional
  approximations}, Graduate Studies in Mathematics, vol.~88, Amer. Math. Soc.,
  Providence, RI, 2008.

\bibitem{BrPeTo:Cuntz}
N.P. Brown, F.~Perera, and A.S. Toms, \emph{The {C}untz semigroup, the
  {E}lliott conjecture, and dimension functions on {$C^*$}-algebras}, J. Reine
  Angew. Math. \textbf{621} (2008), 191--211.

\bibitem{Mitacs2012}
K.~Carlson, E.~Cheung, I.~Farah, A.~Gerhardt-Bourke, B.~Hart, L.~Mezuman,
  N.~Sequeira, and A.~Sherman, \emph{Omitting types and {AF} algebras},
  arXiv:1212.3576, 2012.

\bibitem{CoElIv}
K.~Coward, G.~A. Elliott, and C.~Ivanescu, \emph{The {C}untz semigroup as an
  invariant for {C}$^*$-algebras}, J. reine angew. Math. \textbf{623} (2008),
  161--193.

\bibitem{Dix:Some}
J.~Dixmier, \emph{On some {$C\sp{\ast} $}-algebras considered by {G}limm}, J.
  Functional Analysis \textbf{1} (1967), 182--203.

\bibitem{EllTo:Regularity}
G.A. Elliott and A.S. Toms, \emph{Regularity properties in the classification
  program for separable amenable {$C\sp *$}-algebras}, Bull. Amer. Math. Soc.
  \textbf{45} (2008), no.~2, 229--245.

\bibitem{Fa:Graphs}
I.~Farah, \emph{Graphs and {CCR} algebras}, Indiana Univ. Math. Journal
  \textbf{59} (2010), 1041--1056.

\bibitem{FaHaSh:Model2}
I.~Farah, B.~Hart, and D.~Sherman, \emph{Model theory of operator algebras
  {II}: Model theory}, preprint, arXiv:1004.0741, 2010.

\bibitem{FaHaSh:Model1}
\bysame, \emph{Model theory of operator algebras {I}: {S}tability}, Proc.
  London Math. Soc. (to appear).

\bibitem{FaHaTiKa:Simple}
I.~Farah, D.~Hathaway, A.~Tikuisis, and T.~Katsura, \emph{A simple
  $\mathrm{C}^*$-algebra with finite nuclear dimension which is not {${\mathcal
  Z}$}-stable}, preprint, arXiv:1301.5030, 2013.

\bibitem{FaKa:Nonseparable}
I.~Farah and T.~Katsura, \emph{Nonseparable {UHF} algebras {I}: {D}ixmier's
  problem}, Adv. Math. \textbf{225} (2010), no.~3, 1399--1430.

\bibitem{FaSh:Dichotomy}
I.~Farah and S.~Shelah, \emph{A dichotomy for the number of ultrapowers},
  Journal of Mathematical Logic \textbf{10} (2010), 45Ð81.

\bibitem{FaToTo:Turbulence}
I.~Farah, A.S. Toms, and A.~T\"ornquist, \emph{Turbulence, orbit equivalence,
  and the classification of nuclear {C*}-algebras}, J. Reine Angew. Math. (to
  appear).

\bibitem{Glimm:On}
J.G. Glimm, \emph{On a certain class of operator algebras}, Trans. Amer. Math.
  Soc. \textbf{95} (1960), 318--340.

\bibitem{HaWa}
R.~Haydon and S.~Wassermann, \emph{A commutation result for tensor products of
  {C*}-algebras}, Bull. London Math. Soc. \textbf{5} (1973), 283--28.

\bibitem{Kat:Non-separable}
T.~Katsura, \emph{Non-separable {AF}-algebras}, Operator Algebras: The Abel
  Symposium 2004, Abel Symp., vol.~1, Springer, Berlin, 2006, pp.~165--173.

\bibitem{Ku:Book}
K.~Kunen, \emph{An introduction to independence proofs}, North--Holland, 1980.

\bibitem{Phi:Classification}
N.C. Phillips, \emph{A classification theorem for nuclear purely infinite
  simple {$C\sp *$}-algebras}, Doc. Math. \textbf{5} (2000), 49--114
  (electronic).

\bibitem{Phi:Simple}
\bysame, \emph{A simple separable {$C\sp *$}-algebra not isomorphic to its
  opposite algebra}, Proc. Amer. Math. Soc. \textbf{132} (2004), no.~10,
  2997--3005.

\bibitem{RoLaLa:Introduction}
M.~R\o rdam, F.~Larsen, and N.J. Laustsen, \emph{An introduction to {K}-theory
  for {C}$^*$ algebras}, London Mathematical Society Student Texts, no.~49,
  Cambridge University Press, 2000.

\bibitem{Ror:Classification}
M.~R{\o}rdam and E.~St{\o}rmer, \emph{Classification of nuclear {$C\sp
  *$}-algebras. {E}ntropy in operator algebras}, Encyclopaedia of Mathematical
  Sciences, vol. 126, Springer-Verlag, Berlin, 2002, Operator Algebras and
  Non-commutative Geometry, 7.

\bibitem{sato09b}
R.~Sasyk and A.~T{\"o}rnquist, \emph{Borel reducibility and classification of
  von {N}eumann algebras}, Bulletin of Symbolic Logic \textbf{15} (2009),
  no.~2, 169--183.

\bibitem{Sh:Classification}
Saharon Shelah, \emph{Classification of first order theories which have a
  structure theorem}, Bull. Amer. Math. Soc. (N.S.) \textbf{12} (1985), no.~2,
  227--232.

\bibitem{Sh:E59}
\bysame, \emph{General non-structure theory}, preprint, arXiv:1011.3576, 2010.

\bibitem{ShUs:928}
Saharon Shelah and Alex Usvyatsov, \emph{Unstable classes of metric
  structures}, preprint, available at
  http://ptmat.fc.ul.pt/$\sim$alexus/papers.html.

\bibitem{Wid:Nonisomorphic}
Harold Widom, \emph{Nonisomorphic approximately finite factors}, Proc. Amer.
  Math. Soc. \textbf{8} (1957), 537--540.

\end{thebibliography}
\end{document}